\documentclass{article}
\usepackage{graphicx}
 \usepackage{mathptmx}
\usepackage{amsmath, amstext,amssymb,amsfonts}
\usepackage[english]{babel}
\usepackage{delarray}
\usepackage{mathptmx}
\usepackage{amsmath, amstext,amssymb,amsfonts}
\usepackage[english]{babel}
\usepackage{delarray}

\DeclareFontFamily{U}{mathx}{\hyphenchar\font45}
\DeclareFontShape{U}{mathx}{m}{n}{
      <5> <6> <7> <8> <9> <10>
      <10.95> <12> <14.4> <17.28> <20.74> <24.88>
      mathx10
      }{}
\DeclareSymbolFont{mathx}{U}{mathx}{m}{n}
\DeclareFontSubstitution{U}{mathx}{m}{n}
\DeclareMathAccent{\widecheck}{0}{mathx}{"71}
\DeclareMathAccent{\wideparen}{0}{mathx}{"75}

\numberwithin{equation}{section}
\newtheorem{theorem}{Theorem}[section] \newtheorem{corollary}[theorem]{Corollary}
 \newtheorem{lemma}[theorem]{Lemma}

 \newtheorem{example}[theorem] {Example}
 \newcommand{\bo}{{\cal{B}}(\R)}      \newcommand{\M}{M({\mathbb R})}      \newcommand{\B}{B({\mathbb R})}    
   \newcommand{\R}{\mathbb R}    \newcommand{\Co}{\mathbb C}    \newcommand{\To}{\mathbb T}             \newcommand{\Z}{\mathbb Z} 
   \newcommand{\Bg}{B(G)}   \newcommand{\Mg}{M(G)}   \newcommand{\Mgg }{M(\widehat{G})}  \newcommand{\bog}{{\cal{B}}(G)} \newcommand{\bogh}{{\cal{B}}(\widehat{G})}  \newcommand{\Bgg}{B(\widehat{G})}

\date{}
\title{On the imaginary part of the characteristic function} 
\author{\large{ Saulius  Norvidas}}
\date{\footnotesize Institute of Data Science and Digital Technologies, Vilnius University, \\ Akademijos str. 4, Vilnius LT-04812, Lithuania\\
 ({\rm{e-mail: norvidas{@}gmail.com}})}
\begin{document}

\maketitle
 {{ {\bf Abstract}}}
Suppose that $f$ is the characteristic function of a probability measure on the real line $\R$. In this paper,  we deal with the following   problem posed by N.G. Ushakov: Is  it  true  that  $f$ is  never  determined by  its  imaginary  part $\Im f$? In other words, is  it  true  that for any characteristic function $f$ there exists a characteristic function $g$ such that $\Im f\equiv \Im g$ but $ f\not\equiv  g$? We study this question  in the more general case of the characteristic function defined on an arbitrary  locally compact abelian group.  A characterization of what  characteristic functions  are uniquely determined by their imaginary parts are given. As a consequence of this characterization, we obtain  that several  frequently used characteristic functions  on the classical locally compact abelian groups   are uniquely determined by their imaginary parts.



\section{ Introduction }
{\large{
Let  $\bo$  denote the usual $\sigma$-algebra of all Borel  subsets of the real line $ \R$, and let  $M(\R)$ be  the  Banach algebra of bounded regular complex-valued  Borel measures $\mu$ on $\R$. $M(\R)$ is  equipped with the usual total variation norm $\|\mu\|$.  We define the Fourier (Fourier-Stieltjes) transform of\ $\mu\in\M$  by
\[
\hat{\mu}(\xi)=\int_{\R}e^{-i\xi x} d\mu(x),\quad \xi\in\R.
\]
The family of  these functions  \ $\hat{\mu}$\ forms  the  Fourier-Stieltjes algebra \ $\B$.\  The norm in $\B$  is defined by
\[
\|\hat{\mu}\|_{\B}=\|\mu\|_{\M}.
\]
If $\mu\in\M$ is nonnegative measure and $\|\mu\|=1$, then in the language of probability theory,  $\mu$ and   $f(\xi):= \hat{\mu}(-\xi)$, $\xi\in\R$,  are called a probability measure and its   characteristic function, respectively.

 In this paper,  we deal with several  examples and certain assertions demonstrating the relationship between the characteristic function $f$  and its  imaginary part $\Im f$.
 It is well known that $f$ is not determined by $|f|$.  More precisely,  there exist two different real-valued characteristic functions $f$ and $g$ such that $|f|=|g|$ everywhere (see \cite[ p. 506]{4} and \cite[ p. 265]{11}).
  Next,  any characteristic function $f\not\equiv 1$ is never determined  by its real part $\Re f$, i.e., there exists the characteristic function $g$ such that $\Re f\equiv \Re g$ but $ f\not\equiv  g$ (see \cite[ p. 259]{11}).    In this context is natural to ask whether the same is true for $f$ and $\Im f$. The following question is given in \cite[ p.  334]{11} as an unsolved problem:
\begin{gather}\label{1}
Is \ it \ true \ that \ the \  characteristic \ function \ is \ never \ determined \nonumber\\
\ by \ its \ imaginary \ part?
\end{gather}
In other words, is  it  true  that for any characteristic function $f$, there exists a characteristic function $g$ such that $\Im f\equiv \Im g$ but $ f\not\equiv  g$?
It is not difficult to check that the answer to this question is no. For example, if $\Im f(x)= \Im g(x)=\sin x$,  then necessarily $f(x)=g(x)=e^{ix}$. This can  be verified directly by simple arguments. However, we refer to our theorem 1.

We note that if $f$ is the  characteristic function  of a probability measure $\mu\in\M$, then
\[
\Re f(x) =\frac12 \int_{\R}\bigl(e^{ixt}+e^{-ixt}\bigr)\,d\mu(t)\quad {\text{and}}\quad \Im f(x)=\frac1{2i}\int_{\R}\bigl(e^{-ixt}-e^{ixt}\bigr)\,d\mu(t).
 \]
Therefore,  both $\Re f$ and $\Im f$ are elements of  $\B$. Moreover, $\Re f$ is an even characteristic function  with $\|\Re f\|_{\B}=1$. On the other hand, $\Im f$  is  odd,  $\|\Im f\|_{\B}\le 1 $,   and $\Im f$ is positive definite if and only if $\Im f \equiv 0$.

 We study  the question (1) in the more general case of the characteristic function defined on an arbitrary  locally compact abelian group $G$  (see, e.g., \cite{6}, \cite{10}, and the next section of this paper  for the exact definitions of algebras $\bog$, $\Mg$,  $\Bg$, and other  background of harmonic analysis on groups).

  Let us start with a  question: which a function  $\varphi:G\to \R$  might serve as the imaginary part of the characteristic function $f$ on $G$? Keeping in mind what was said above about the imaginary part of the characteristic function, the following theorem gives a complete description of these functions $\varphi$.

\begin{theorem}\label{1}
Assume that  $\varphi\in\Bg$ is real-valued and   odd.  If   $\|\varphi\|_{\Bg}\le 1$, then there exists a characteristic function $f$  on $G$  such that  $\Im f\equiv \varphi$.
\end{theorem}

As a consequence,  we obtain  a characterization of what  characteristic functions  are uniquely determined by their imaginary parts.
\begin{theorem}\label{2}
Let $f:G\to \Co$ be a characteristic function. Then  $\Im f$ completely determines $f$ if and only if $\|\Im f\|_{\Bg}=1$.
\end{theorem}
The following  corollary can be used to construct characteristic functions  which are uniquely determined by their imaginary parts.
 We denote by $\widehat{G}$ the dual group of $G$ (see the  next section   for the exact definition). For  $U\subset \widehat{G}$, we write   $-U=\{t\in G: \ -t\in U\}$.
\begin{corollary}\label{3}
Suppose that  $f:G\to \Co$ is the characteristic function of a probability measure  $\mu\in M(\widehat{G})$.
Let $U\in\bogh$, and assume that $U\cap(-U)=\emptyset$.  If
\begin{equation}\label{1.2}ƒ
\mu(U)=1,
\end{equation}
 then  $f$ is   completely determined by $\Im f$ .
\end{corollary}
Now it is easy to verify   that several frequently used characteristic functions  on the classical groups $\R$ , $\Z$, $\To$, and $\R^{n}$ are uniquely determined by their imaginary parts  (see, e.g.,   a list of several probability distributions and  their characteristic functions    in \cite[ p.p. 282-329]{11}). At the end of this chapter, we recall  some distributions and  their characteristic functions of this type.

\begin{example}\label{4}ƒ
The  group $G=\R$.

(A) \ The characteristic  functions  of the following probability distributions  are uniquely determined by their imaginary parts:

the arcsine distribution,  the Bessel distribution, the beta distribution, the gamma distribution, the  hyperexponential distribution, the standard Levy distribution, the Maxwell distribution, the Pareto distribution, and the $\chi^2$-distribution  with  density function
\[
p(t)=\begin{array}\{{rl}.
\frac{1}{2^{n/2}\Gamma(n/2)}t^{(n/2)-1}e^{-t/2}   , &  t>0, \\
0, & \mbox{otherwise,}
\end{array}
\]
where $n$  is a positive  integer.

(B) \ The characteristic  functions  of the following probability distributions are not uniquely determined by their imaginary parts:

 the normal distribution,  the Laplace distribution,  and the Cauchy distribution.

(C) \ The characteristic  functions $f$ of the uniform distribution and the triangular (Simpson) distribution with density  defined by
\[
p(t)=\begin{array}\{{rl}.
\frac1{b-a}, & a\le t\le b, \ a<b, \\
0, & \mbox{otherwise,}
\end{array}
\]
and by
\[
p(t)= \begin{array}\{{rl}.
\frac4{(a-b)^2}\Bigl( \frac{b-a}2-|t-\frac{a+b}2|\Bigr), & a\le t\le b, \ a<b,\\
0, & \mbox{otherwise,}
\end{array}
\]
respectively, are uniquely determined by $\Im f$ if and only if \ $0\not\in [a, b]$.
\end{example}
\begin{example}\label{5}ƒ
The  groups $G=\R^n$, $n>1$.

The  characteristic functions  of the  multivariate distributions on $\R^n$ inherit  in most cases   the same properties as the appropriate univariates distributions.   For example, the characteristic function of the Dirichlet distribution on $\R^n$, i.e., the multivariate generalization of the beta distribution (see, e.g., \cite[ Chapter 49]{7}), and multivariate Pareto distributions of the first and the second kinds (see, e.g., \cite[ Chapter 6]{1}) are uniquely determined by their imaginary parts.
 On the other hand, the characteristic functions of the normal distribution,  the Laplace distribution,  and the Cauchy distribution  are not uniquely determined by their imaginary parts.
\end{example}

\begin{example}\label{6}ƒ
The  group $G=\To$.

In these examples we  consider the circle group $\To$ as     the interval $\To=\R/(2\pi \Z)=[0,2\pi)$  with addition mod $2\pi$ and with the usual topology inherit from $\R$.  Then $\widehat{\To}=\Z$. The Bochner theorem [6, p. 293]  states now that a function $f: \To\to \Co$ is the characteristic function of a probability measure on $\Z$ if and only if
\begin{equation}\label{3}ƒ
f(x)=\sum_{k\in\Z}\alpha_ke^{ikx},
\end{equation}
where $\alpha_k\ge 0$ and $\sum_{k} \alpha_k =1$.  Let $A$ denote the set of all $k\in\Z$ in (3) such that $\alpha_k> 0$. Then (3) is uniquely determined by its imaginary part if and only if  $A\cap(-A)=\emptyset$. Therefore, the characteristic functions of the negative binomial distribution, the Poisson distribution, the binomial distribution, and the hypergeometric distribution are not uniquely determined by their imaginary parts.
\end{example}

\begin{example}\label{7}ƒ
The  group $G=\Z$.

In this case,  the set of characteristic functions on $\Z$ coincides with the family of positive definite sequences $\{\xi_k: \ \xi_k\in\Z\}$ of complex numbers $\xi_k$ such that $\xi_0=1$. Since $\Z=\widehat{\To}$, we see that any such a sequence is the Fourier coefficients sequence of a probability distribution on $\To$. In other words, this $\{\xi_k: \ \xi_k\in\Z\}$ coincides with the Fourier coefficients sequence of a periodic distribution on $\R$ with period $2\pi$. A periodic distribution can be obtained by the periodization ( by the "wrapping" ) around the unit circle $\To=\R/(2\pi \Z)$ of a distribution on $\R$ (see, e.g., \cite[ Chapter 3]{9}). For example, the characteristic functions $f$ of the wrapped Cauchy distribution, the wrapped normal distribution, and he wrapped exponential distribution  are not uniquely determined by $\Im f$. On the other hands, according to the corollary 1.3, it is easy to create the probability distributions on $\To$ such that their characteristic functions, i.e., the characteristic sequences are uniquely determined by their imaginary parts.
\end{example}

\section{PRELIMINARIES}
\label{sec:1}

We recall and  introduce some terminology and notation. Throughout $G$ will be a locally compact abelian group with  dual group $\widehat{G}$  consisting of the continuous homomorphisms (or characters) $\gamma: G \to \To$, where $\To$  is the multiplicative unit circle group in $\Co$. The dual $\widehat{G}$  is also a locally compact abelian group. We write the group operations in $G$ and $\widehat{G}$  additively. Let us denote the identity elements of G and $\widehat{G}$ by $0$, and write the action of a character $\gamma\in\widehat{G}$ at the point $x\in G$ as $\gamma(x)=(x, \gamma)$. Note that $(-x, \gamma)=(x, -\gamma)=\overline{(x, \gamma)}$.

Let  $\bog$ denote the smallest $\sigma$-algebra that contains all closed subsets of $G$. The elements  of $\bog$ are called the Borel sets of $G$. A measure $\mu$ on $G$ is a complex-valued function on $\bog$, which is $\sigma$-additive set function , i.e.,
 \[
 \sum_{k\in \Z}\mu(E_k)=\mu(\sum_{k\in \Z}E_k)
 \]
  if $\{E_k\}$ are pairwise disjoint elements of $\bog$, and $\mu(K)$ is finite for all compact subset $K$ of $G$. If $\mu(E)$ is finite for all $E\in \bog$, then $\mu$ is called a bounded measure.  Let $\Mg$ denote the Banach algebra of all complex-valued, regular and bounded   measures on $\bog$. A norm is introduced in M(G) by defining
\[
\|\mu\|=  \sup \sum_k|\mu(E_k) |,
\]
where the sup being taken over all finite collections of pairwise disjoint Borel sets $E_k$ such that $\cup_k E_k=G$.

Let $m$ be a nontrivial Haar measure defined on $G$. Then $L^p(G)$, $1\le p\le \infty$,   denotes the  usual Lebesgue  space with respect to $m$ of  $\bog$-measurable, complex valued functions defined on G. There corresponds to each $\varphi\in L^1(G)$  an unique measure $\mu_{\varphi}\in \Mg$  defined by
\[
\mu_{\varphi}(E) = \int_E\varphi(x)\,dx,
\]
  where  $ E\in\bog$. By the Radon-Nikodym theorem $\mu=\mu_{\varphi}$ for some $\varphi\in  L^1(G)$ if and only if $\mu$ is absolutely continuous with respect to $m$. Recall that if $\mu$ is a probability measure and there exists $\varphi\in L^1(G)$ such that  $\mu=\mu_{\varphi}$, then $\varphi$ is called the density (density function) of $\mu$  (with respect to $m$). Thus,   $L^1(G)$ can be identified with the closed ideal in $\Mg$ of all measures absolutely continuous  with respect to the Haar   measure $m$ on $G$.

For $\mu\in\Mg$,  Fourier-Stieltjes transform of $\mu$  is the function $\widehat{\mu}$ defined on $\widehat{G}$  by
\[
\widehat{\mu}(\gamma)=\int_G \overline{(x, \gamma)}\,d\mu(x)= \int_G (-x, \gamma)\,d\mu(x),
\]
$\gamma\in\widehat{G}$. By the Pontryagin duality theorem that every locally compact abelian group is (isomorphic to) the dual group of its dual group, i.e., $\widehat{(\widehat{G})}=G$, we will also write
\[
\widehat{\omega}(x)=\int_{\widehat{G}} \overline{(x, \gamma)}\,d\omega(\gamma)= \int_{\widehat{G}} \overline{(\gamma,x)}\,d\omega(\gamma),\quad x\in G,
\]
for any $\omega\in \Mgg$. If $\varphi\in L^1(G)$, then
\[
\widehat{\varphi}(\gamma)=\int_G \overline{(x, \gamma)}\varphi(x)\,dm(x),\quad \gamma\in\widehat{G},
\]
is said to be the Fourier transform of $\varphi$.

The set of Fourier-Stieltjes transforms $\widehat{\mu}$,  $\mu\in \Mg$,  is a function algebra $\Bgg$ on $\widehat{G}$, the  Fourier-Stieltjes algebra of $\widehat{G}$, with the ordinary pointwise algebraic operations.  The norm in $\Bgg$  is defined by $\|\hat{\mu}\|_{\Bgg}=\|\mu\|_{\Mg}$. The closed ideal\ $A(\widehat{G})=\{\hat{\varphi}: \ \varphi\in L^1(G)\}$ in $\Bgg$ is called the Fourier algebra of $\widehat{G}$.

A function $f:G\to\Co$ is said to be positive definite if
\[
\sum_{j, k=1}^n f(x_j-x_k)c_j{\overline{c}}_k\ge 0
\]
holds for all finite sets of complex numbers $c_1,\dots, c_n$ and points $x_1,\dots,x_n\in G$.     The Bochner theorem (see, e.g.,  \cite[ p. 121]{3}, \cite[ p. 293]{6}, and \cite[ p. 71]{8})  states that a continuous function $f: G\to\Co$  is positive definite if and only if there exists a nonnegative $\mu\in \Mgg$ such that $f= \hat{\mu}$. If, in addition, in this theorem we have that $\|\mu\|=1$, then  $f(0)=1$ and $f$ is called the    characteristic function of $\mu$. Let $P(G)$ denote the set of all continuous, positive definite functions defined on $G$.

 We call a measure $\mu\in\Mg$ symmetric  if $\mu(-A)=\mu(A)$ for all $A\in\bog$. Of course,  if $\mu$ is a probability measure on $G=\R$, then, in terms of the probability theory, $\mu$ is symmetric iff the distribution function corresponding to $\mu$ is symmetric (see \cite[ p. 30]{8}). We say that $\mu\in\Mg$ antisymmetric  if $\mu(-A)=-\mu(A)$,  $A\in\bog$. Given $\mu\in\Mg$, we associate with $\mu$ two other measures, defined  by
\[
\mu_s(A)= \frac12\Bigl(\mu(A)+\mu(-A)\Bigr), \qquad \mu_a(A)= \frac12\Bigl(\mu(A)-\mu(-A)\Bigr)
\]
for   $A\in\bog$. We say that   $\mu_s$ and $\mu_a$  are the   symmetric part and antisymmetric  part of $\mu$, respectively.  It is obvious that $\mu_s$ and $\mu_a$ are symmetric and antisymmetric measures in $\Mg$, respectively.

Let $\mu\in \Mg$. A set $A\in\bog$ is  positive (resp. negative, null) set for $\mu$ if $\mu(E)\ge 0$ (resp. $\mu (E)\le 0 $, $\mu(E)=0$) for all $E\in\bog$ such that $E\subset A$ (see, e.g., \cite[ p. 86]{5}). The measures $\mu$ and $\eta$ are called mutually singular (or  $\mu$  is singular with respect  $\eta$, or vice versa) if there exist $A, B\in\bog$ such that
\begin{equation}\label{4}ƒ
A\cap B=\emptyset, \qquad A\cup B=G,
\end{equation}
  $A$ is null set for $\mu$, and $B$ is null set for $\eta$.

We call any real-valued measure $\mu$ in $\Mg$  a signed measure. For every signed measure $\mu\in \Mg$,ƒ  the following decomposition is true (see, e.g., \cite[ p.p. 175-176]{2} and \cite[ p.p. 86-87]{5}).
\begin{theorem}\label{8}(The Hahn-Jordan decomposition theorem.)
If $\mu\in \Mg$ is a signed measure, then there exist unique  mutually singular positive measures $\mu^+$ and $\mu^-$ in $\Mg$ such that $\mu=\mu^+-\mu^-$ (the Jordan decomposition). In addition, there exist  $A^+, A^-\in\bog$ such that $A^+\cap A^-=\emptyset, \qquad A^+\cup A^-=G$, measures $\mu^+$ and  $\mu^-$ are supported on $A^+$ and $A^-$, respectively, in the sense that  $\mu^+(A^-)=\mu^-(A^+)=0$ (the Hahn decomposition). If $A^+_1$ and $A^-_1$  is another such pair, then $A^+\triangle A^+_1$ and $A^-\triangle A^-_1$ are null sets for $\mu$.
For all $E\in\bog$,
\begin{equation}\label{5}
\mu^+(E)=\mu(E\cap A^+) \qquad and \qquad \mu^-(E)=-\mu(E\cap A^-) .
\end{equation}
\end{theorem}
ƒ
In the case  of symmetric or antisymmetric measures, this Hahn-Jordan decomposition can be strengthened. We  need one of that  statement. For completeness, we give also its proof.

\begin{lemma}\label{9}
Let  $\mu\in \Mg$ be a signed measure,  let  $\mu=\mu^+-\mu^-$ be its Jordan decomposition, and suppose that $A^+$ and $A^-$ form its Hahn decomposition. Set $V= A^+\cap(-A^-)$. If $\mu$ is, in addition,  antisymmetric, then
\begin{equation}\label{6}ƒ
V\cap(-V)=\emptyset,
\end{equation}
\begin{gather}\label{7}ƒ
\mu^+(E)=\mu(E\cap A^+)=\mu(E\cap V), \quad \mu^-(E)=-\mu(E\cap A^-)=-\mu(E\cap (-V))
\end{gather}
for all $E\in\bog$. Moreover,
\begin{equation}\label{8}ƒ
\mu^+(V)=\|\mu^+\|=\|\mu^-\|=\mu^-(-V)=\frac12\|\mu\|.
\end{equation}
\end{lemma}

{\it{ Proof.}} \

The Hahn decomposition of $\mu$ shows that $A^+$ and $A^-$ are  positive and negative sets for $\mu$, respectively. Since $\mu$ is antisymmetric, we have that $-A^+$ and $-A^-$  are  negative and positive for $\mu$, respectively. Moreover, is clear that
\[
(-A^+)\cap(- A^-)= - (A^+\cap A^-)=\emptyset \quad {\text{and}}\quad (-A^+)\cup(- A^-)=-(A^+\cup A^-)= G.
 \]
 Therefore,  $-A^-$ and $-A^+$ is also a Hahn decomposition for $\mu$. Thus, by theorem 3, we see that   $A^+\triangle (-A^-)$ and $(-A^+)\triangle A^-$ both are null sets for $\mu$.

 Set
\begin{equation}\label{9}ƒ
T_1= A^+\setminus\Bigl(A^+\triangle(-A^-)\Bigr), \qquad T_2= A^-\setminus\Bigl(A^-\triangle(-A^+)\Bigr).
\end{equation}
Then
\begin{equation}\label{10}ƒ
T_1\cap T_2=\emptyset
\end{equation}
  and
\begin{gather}\label{11}ƒ
\mu(T_1)=\mu(A^+)=\mu^+(A^+)=\mu^+(G)=\|\mu^+\|, \nonumber\\ \mu(T_2)=\mu(A^-)=-\mu^-(A^-)=-\mu^-(G)=\|\mu^-\|.
\end{gather}
Now a direct computation shows that
\begin{equation}\label{12}ƒ
T_1= A^+\cap(-A^-)=V,\quad  {\text{and}} \quad T_2=(- A^+)\cap(A^-)=-\Bigl( A^+\cap(-A^-)\Bigr)=-V.
\end{equation}
Therefore, (10) implies (6). By combining (11) with (12) and using the fact that $\mu$ is antisymmetric, we obtain (8). According to (9) and by the above remark that $A^+\triangle (-A^-)$ and $(-A^+)\triangle A^-$ both are null sets for $\mu$, we have
\[
\mu(E\cap A^+)=\mu(E\cap T_1)\quad {\text{and}} \quad \mu(E\cap A^-)=\mu(E\cap T_2 )
\]
for all $E\in\bog$. Combining these with (5) and (12), we obtain (7).  Lemma 1 is proved.
\section{PROOFS}
\label{sec:2}
{\it{ Proof of  Theorem 1.}} \

 By the definition of $\Bg$, there exists a $\mu\in \Mgg$ such that $\varphi=\widehat{\mu}$. Let
 \begin{equation}\label{13}ƒ
 \mu=\nu+i \eta,
 \end{equation}
  where $\nu$ and $\eta$ are signed measures equal to $\Re \mu$ and $\Im \mu$, respectively. Since $\varphi$ is real-valued, we have
 \begin{gather}\label{14}ƒ
\varphi(x)=\Re \widehat{\mu}(x)=\Re\biggl(\int_{\widehat{G}}\overline{(x, \gamma)}\,d\mu(\gamma)\biggr)\nonumber\\
= \int_{\widehat{G}}\Re\Bigl(\overline{(x, \gamma)}\Bigr)\,d\nu(\gamma) - \int_{\widehat{G}}\Im\Bigl(\overline{(x, \gamma)}\Bigr)\,d\eta(\gamma).
\end{gather}
For $x\in G$, the functions
\[
x\to \Re\Bigl(\overline{(x, \gamma)}\Bigr)=\frac12\Bigl(\overline{(x, \gamma)}+(x, \gamma)\Bigr)=\frac12\Bigl((-x, \gamma)+(x, \gamma)\Bigr)
\]
and
\[
x\to \Im\Bigl(\overline{(x, \gamma)}\Bigr)=\frac1{2i}\Bigl(\overline{(x, \gamma)}-(x, \gamma)\Bigr)=\frac1{2i}\Bigl((-x, \gamma)-(x, \gamma)\Bigr)
\]
are even and odd, respectively. Therefore, using the fact that $\varphi$ is odd, we obtain from (14) that
 \begin{equation}\label{15}ƒ
\varphi(x)= \frac1{2i}\int_{\widehat{G}}\Bigl((-x, \gamma)-(x, \gamma)\Bigr)\,d\eta(\gamma).
\end{equation}
On the other hand, if $\eta_a$ denotes the the antisymmetric part of $\eta$, then
\begin{gather}
\widehat{\eta_a}(x)= \int_{\widehat{G}}\overline{(x, \gamma)}\,d\eta_a(\gamma)=\frac12\biggl(\int_{\widehat{G}}\overline{(x, \gamma)}\,d\eta(\gamma)-\int_{\widehat{G}}\overline{(-x, \gamma)}\,d\eta(\gamma)\biggr)\nonumber\\
=\frac12\int_{\widehat{G}}\Bigr((-x, \gamma)-(x, \gamma)\Bigr)\,d\eta(\gamma)\nonumber.
\end{gather}
Combining this with (14) and (15), we get
 \begin{equation}\label{16}ƒ
i\varphi(x)=\widehat{\eta_a}(x)=\int_{\widehat{G}}\overline{(x, \gamma)}\,d\eta_a(\gamma),
\end{equation}
$x\in G$.  Therefore, if $\eta_a=\eta_a^+-\eta_a^-$ is the Jordan decomposition of $\eta_a$, then
\begin{equation}\label{17}ƒ
i\varphi=\widehat{\eta_a^+}-\widehat{\eta_a^-}.
\end{equation}
Set
\[
|\eta_a|=\eta_a^++\eta_a^-.
\]
  Since $\eta_a$ is asymmetric, it follows that the measure $|\eta_a|$ is symmetric. Indeed, if $E\in \bogh$, then according to lemma 1 and (7), we get
\begin{gather}
 |\eta_a|(-E)= \eta_a^+(-E)+\eta_a^-(-E)= \eta_a(-E\cap V)-\eta_a(-E\cap(-V))= \nonumber\\
 -\eta_a(-( E\cap V))) +\eta_a(-(-E\cap(-V)))= - \eta_a(E\cap (-V)) +\eta_a(E\cap V)= \nonumber\\
 \eta_a^-(E)+\eta_a^+(E)=  |\eta_a|(E),\nonumber
\end{gather}
where the set  $V\in \bogh$  is defined by lemma 1 for asymmetric measure $\eta_a$ on $\widehat{G}$ (instead of $\mu$ on $G$ in lemma 1). Let $\omega$ denote the Fourier-Stieltjes transform of $|\eta_a|$. By Bochner's  theorem,  $\omega$ is  continuous real-valued and positive definite. Therefore,  $\omega$ is also even function. Moreover, (16) and (17) imply that
\begin{equation}\label{18}ƒ
\omega+i\varphi = \widehat{|\eta_a|}+ \widehat{\eta_a}=2\widehat{ (\eta_a^+)}.
\end{equation}
 Hence   $\omega+i\varphi\in P(G)$. Now using lemma 1, combining  (8) (in the case where $\mu=\eta_a$)  with (16) and (18), we get
\[
(\omega+i\varphi)(0)=2\widehat{(\eta_a^+)}(0)=2\|\eta_a^+\|=\|\eta_a\|=\|\varphi\|_{\Bg}.
\]
The requirement of our theorem on $\varphi$ is that $\|\varphi\|_{\Bg}\le 1$. Therefore, if we  take any  real-valued even  $\sigma\in P(G)$ such that $\sigma(0)=1-\|\varphi\|_{\B}$, then
\begin{equation}\label{19}ƒ
\sigma+ \omega+i\varphi
\end{equation}
is the characteristic function with  the prescribed imaginary part $\varphi$. Theorem 1 is proved.

{\it{ Proof of  Theorem 2.}} \

Let $f=\psi+i\varphi$, where $\psi =\Re f$ and $\varphi= \Im f$. Assume first that  $\|\varphi\|_{\B}<1$. Then, as  in the final part of the proof  of  theorem 1, we see  that there exist infinitely many disjoint  characteristic functions (19) with   a prescribed imaginary part $\varphi$.

 Conversely, suppose now that $\|\varphi\|_{\B}=1$. Since $\varphi$ is real-valued and odd, we obtain, as  in the  proof  of  theorem 1,  that there exists an antisymmetric measure $\eta_a$ on $\widehat{G}$ with the Jordan decomposition  $\eta_a=\eta_a^+-\eta_a^-$ such that   (17) is satisfied. Suppose that  $V$ is defined by lemma 1 for  $\eta_a$ instead of $\mu$. Then we get from (8) that
\begin{equation}\label{20}ƒ
\eta_a^+(V)=\|\eta_a^+\|=\|\eta_a^-\|=\eta_a^-(-V)=\frac12\|\eta_a\|=\frac12\|\varphi\|_{\Bg}=\frac12.
\end{equation}
The real part $\psi=\Re f$ is also the characteristic function. Then there exist  a symmetric probability measure $\tau$ on $\widehat{G}$ such that $\widehat{\tau}=\psi$. Now set
\begin{equation}\label{21}ƒ
\mu=\tau+\eta_a=\tau+\eta_a^+-\eta_a^-.
\end{equation}
 We conclude from (16) that $f=\widehat{\mu}$. Therefore, $\mu$ is a probability measure. Using (20) and the fact that $\eta_a^+(-V)=0$, we get
 \[
 0\le \mu(-V)=\tau(-V)-\eta_a^-(-V)=\tau(-V)-1/2.
  \]
Hence   $\tau(-V)\ge 1/2$.   Therefore,
\begin{equation}\label{22}ƒ
\tau(V)=\tau(-V)=\frac12 \qquad {\text{and}} \qquad \tau\Bigl(\R\setminus(V\cap(-V)\Bigr)=0,
\end{equation}
since $\tau$ is a symmetric probability measure. Define  the following two non-negative measures $\tau_1$ and $\tau_2$ on $\widehat{G}$ by
\begin{equation}\label{23}ƒ
\tau_1(E)=\tau(E\cap V) \quad {\text{and}}\quad \tau_2(E)=\tau(E\cap(-V)),
\end{equation}
$E\in\bogh$. Now (22) implies that $\tau_1$ and $\tau_2$ are mutually singular. Furthermore,
\begin{equation}\label{24}
\tau=\tau_1+\tau_2.
\end{equation}
 We claim that
\begin{equation}\label{25}
ƒ\tau_1=\eta_a^+  \qquad {\text{and}} \qquad \tau_2=\eta_a^-.
 \end{equation}
At first, we recall that if $A^+$ and $A^-$ is a Hahn decomposition for $\eta_a$, then   $\eta_a^{\pm}$ are defined as $\eta_a^{+}(E)=\eta_a(E\cap A^+)$  and $\eta_a^{-}(E)=-\eta_a(E\cap A^-)$, $E\in \bogh$. Moreover, by lemma 1, we have
\begin{equation}\label{26}
\eta_a^{+}(E)=\eta_a(E\cap V) \qquad {\text{  and}} \qquad \eta_a^{-}(E)=-\eta_a(E\cap (-V))
 \end{equation}
for all $E\in\bogh$. Combining (23) with (26),  we get
\[
0\le \mu(E\cap(- V))=\tau(E\cap (- V))+\eta_a(E\cap (- V))= \tau_2(E)-\eta_a^-(E).
\]
Therefore,  for all $E\in \bogh$,
\begin{equation}\label{27}
 \tau_2(E)\ge \eta_a^-(E).
\end{equation}
From (20) it follows that   $\|\tau_2\|\ge \|\eta_a^-\|=1/2$. On the other hand, (22) and (23) give that $\tau_1(V)=\tau(V)=1/2$. Hence $\|\tau_1\|\ge 1/2$.  Finally, since  $\tau=\tau_1+\tau_2$, where $\tau_1$ and $\tau_2$ are mutually singular, we see that $\|\tau_1\|=\|\tau_2\|=1/2$. Therefore,  $\tau$ is a probability measure. Moreover, (20) and (27) give
\begin{equation}\label{28}
\tau_2=\eta_a^-.
\end{equation}
Next, in light of (20), (26),  (28), and using the facts that $\tau$ and $\eta_a$ are  symmetric and antisymmetric, respectively, we have
\begin{gather}
\tau_1(E)=\tau(E\cap V)=\tau(-(E\cap V))=\tau((-E)\cap (-V))=\tau_2(-E)=\eta_a^-(-E)\nonumber\\
=-\eta_a(-(E\cap V))=\eta_a(E\cap V)=\eta_a^+(E)\nonumber
\end{gather}
for each $E\in \bogh$. This, together with (28) proves   our claims (25).
Finally, (25) shows that the real part $\psi=\widehat{\tau}$  of $f$   is completely determined by $\eta_a$. Since $i\Im f=\widehat{\eta_a}$, this means that theorem 2 is proved.

{\it{ Proof of Corollary 1.}} \

 By the definition of asymmetric part of $\mu$, we have
\[
\mu_a(U)=\frac12\Bigr(\mu(U)-\mu(-U)\Bigr)=\frac12
\]
 and
 \[
 \mu_a(-U)=\frac12\Bigr(\mu(-U)-\mu(U)\Bigr)=-\frac12.
\]
Hence  $\|\mu_a\|_{\Mgg}=1$. On the other hand, it is easy to see that
\begin{equation}\label{29}
\widehat{\mu_a}= i \Im \widehat{\mu}.
\end{equation}
 Therefore, $\|\Im f\|_{\Bgg}=\|\Im \widehat{\mu}\|_{\Bgg}=\|\widehat{\mu_a}\|_{\Bgg}=\|\mu_a\|_{\Mgg}=1$.  Finally, using theorem 2,  we obtain the statement of corollary 1.

{\it{ Proof of Examples 1 -- 4.}} \

Let $f$ be an arbitrary  characteristic function  mentioned in these examples and such that we claim that $\Im f$ uniquely determine $f$. Then $f$ satisfies the requirements of corollary 1 (for example, if $G=\R^m$, $m=1,2,\dots$,  then  the statement of corollary 1 is satisfied for  $U=\R^m_{+}=(0;\infty)^m$). In the case of all other type of the characteristic functions $f$ from our examples, there exist $E\in\bogh$ such that $E=-E$ and $\mu(E)>0$,  where $\mu\in\Mgg$ is a probability measure such that $\widehat{\mu}=f$. Hence  $\mu_a(E)=0$. Combining this  with  (29), we get
\begin{gather}
\| \Im f\|_{\Bgg}=\|\Im\widehat{\mu}\|_{\Bgg}=\|\widehat{\mu_a}\|_{\Bgg}=\|\mu_a\|_{\Mgg}=|\mu_a(G)|= |\mu_a(G\setminus E)|+|\mu_a(E)| \nonumber\\
= |\mu_a(G\setminus E)|\le   |\mu(G\setminus E)|<1.                              \nonumber
\end{gather}
 Finally, theorem 2   shows that  $f$ of this type are not uniquely determined by $\Im f$.

.



}}
\end{document}